# Представление функции Миттаг-Леффлера через экспоненциальную функцию в случае рациональных производных

**Фикрет А. Алиев, Н.А. Алиев, Н.А. Сафарова**

*Институт Прикладной Математики БГУ, Баку, Азербайджан*
*e-mail: f_aliev@yahoo.com, nihan.aliev@gmail.com  narchis2003@yahoo.com*

**Абстракт**. В данной работе приводится функция Миттаг-Леффлера через экспоненциальную функцию для любых рациональных производных порядка $\frac{m}{n}$, где $m<n$, $n>1$ являются натуральными несократимыми числами (если $n=1$ то и $m$ тоже равняется единице). В отличие от предыдущих работ приведенные формулы не содержат интегралов.

1. **Введение**

В последнее время, много внимания уделяется построению решения задачи Коши и граничных задач для дифференциальных уравнений дробного порядка [1-4]. Полученные результаты в этой области носят чисто теоретический характер, т.е. или доказывается существование или единственность решений соответствующих задач [5-7].

В работе [8,9] приводится решение задачи Коши для случая любых производных вещественного порядка (стационарные линейные системы), которые можно достаточно хорошо приближать к рациональным дробям, числитель и знаменатель, которых являются нечетными натуральными числами. А это создает некоторые погрешности для решения соответствующих дифференциальный уравнений дробного порядка.

Поэтому, в данной работе показывается, что функцию Миттаг-Леффлера можно представить через экспоненциальную функцию для общего случая, т.е., когда порядок производной является любым правильным рациональным числом. Такая схема позволяет исследовать решение задачи для дифференциальных уравнений дробного порядка, шаг (производных) которого меньше единицы.

2. **Некоторые преобразования для функции Миттаг-Леффлера со сдвигом**

Как известно [1], функцию Миттаг-Леффлер со сдвигом можно представить в следующем виде

$$h_{\frac{1}{n}}(x,\rho) = \sum_{k=0}^{\infty} \rho^k \frac{x^{-1+\frac{k+1}{n}}}{\left(-1+\frac{k+1}{n}\right)!}, \quad x \geq x_0 > 0, \qquad (1)$$

где $n$ - натуральное число, $x$ - вещественный аргумент, $n \in N,\ x \in R,\ x \geq x_0 > 0,\ \rho \in C$ - параметр.

Из (1) легко видно, что

$$D^{\frac{1}{n}} h_{\frac{1}{n}}(x,\rho) = \rho h_{\frac{1}{n}}(x,\rho). \qquad (2)$$

Теперь представим функцию (1) в следующем виде

$$\begin{aligned}
h_{\frac{1}{n}}(x,\rho) &= \frac{x^{-1+\frac{1}{n}}}{\left(-1+\frac{1}{n}\right)!} + \rho \frac{x^{-1+\frac{2}{n}}}{\left(-1+\frac{2}{n}\right)!} + \cdots + \rho^{n-1} \frac{x^0}{0!} + \\
&+ \rho^n \frac{x^{\frac{1}{n}}}{\left(\frac{1}{n}\right)!} + \rho^{n+1} \frac{x^{\frac{2}{n}}}{\left(\frac{2}{n}\right)!} + \cdots + \rho^{2n-1} \frac{x}{1!} + \\
&+ \rho^{2n} \frac{x^{1+\frac{1}{n}}}{\left(1+\frac{1}{n}\right)!} + \rho^{2n+1} \frac{x^{1+\frac{2}{n}}}{\left(1+\frac{2}{n}\right)!} + \cdots + \rho^{3n-1} \frac{x^2}{1!} + \cdots + \\
&+ \rho^{sn} \frac{x^{s-1+\frac{1}{n}}}{\left(s-1+\frac{1}{n}\right)!} + \rho^{sn+1} \frac{x^{s-1+\frac{2}{n}}}{\left(s-1+\frac{2}{n}\right)!} + \cdots + \rho^{(s+1)n-1} \frac{x^s}{s!} + \\
&+ \ldots \quad \ldots \ldots \quad \ldots \ldots \quad \ldots \ldots \quad \ldots \ldots \quad \ldots
\end{aligned} \qquad (3)$$

Суммируя (3) по соответствующим столбцам получим:

$$\begin{aligned}
h_{\frac{1}{n}}(x,\rho) &= \left[ \frac{x^{-1+\frac{1}{n}}}{\left(-1+\frac{1}{n}\right)!} + \rho^n \frac{x^{\frac{1}{n}}}{\left(\frac{1}{n}\right)!} + \cdots + \rho^{sn} \frac{x^{s-1+\frac{1}{n}}}{\left(s-1+\frac{1}{n}\right)!} + \cdots \right] + \\
&+ \rho \left[ \frac{x^{-1+\frac{2}{n}}}{\left(-1+\frac{2}{n}\right)!} + \rho^n \frac{x^{\frac{2}{n}}}{\left(\frac{2}{n}\right)!} + \cdots + \rho^{sn} \frac{x^{s-1+\frac{2}{n}}}{\left(s-1+\frac{2}{n}\right)!} + \cdots \right] + \cdots + \\
&+ \rho^{n-1} \left[ \frac{x^0}{0!} + \rho^n \frac{x}{1!} + \cdots + \rho^{sn} \frac{x^s}{s!} + \cdots \right]
\end{aligned} \qquad (4)$$

или

$$h_{\frac{1}{n}}(x,\rho) = \sum_{k=0}^{\infty} \rho^{kn} \frac{x^{-1+k+\frac{1}{n}}}{\left(-1+k+\frac{1}{n}\right)!} + \rho \sum_{k=0}^{\infty} \rho^{kn} \frac{x^{-1+k+\frac{2}{n}}}{\left(-1+k+\frac{2}{n}\right)!} + \cdots + \rho^{n-1} \sum_{k=0}^{\infty} \rho^{kn} \frac{x^k}{k!}.$$

Или же,

$$h_{\frac{1}{n}}(x,\rho) = J_0(x,\rho^n) + \rho J_1(x,\rho^n) + \cdots + \rho^{n-1} J_{n-1}(x,\rho^n), \qquad (5)$$

где

$$J_s(x,\rho^n) = \sum_{k=0}^{\infty} \rho^{kn} \frac{x^{-1+k+\frac{s+1}{n}}}{\left(-1+k+\frac{s+1}{n}\right)!}, \quad s = \overline{0, n-1}. \qquad (6)$$

Формулы (4) и (5) дают возможность связать функцию Миттаг-Леффлера со сдвигом через экспоненциальную функцию. Для этого достаточно свести функции $J_s(x,\rho^n)$ из (6) к экспоненциальной функции.

## 3. Сведение функции $J_s(x,\rho^n)$ к экспоненциальной функции

Сейчас рассмотрим дробный производный порядок $\frac{s+1}{n}$ от функции $J_s(x,\rho^n)$ из (6).

$$\begin{aligned}
D^{\frac{s+1}{n}} J_s(x,\rho^n) &= \sum_{k=0}^{\infty} \rho^{kn} \frac{x^{k-1}}{(k-1)!} = \frac{x^{-1}}{(-1)!} + \rho^n \frac{x^0}{0!} + \rho^{2n} \frac{x}{1!} + \rho^{3n} \frac{x^2}{2!} + \cdots = \\
&= \rho^n \left( \frac{x^0}{0!} + \rho^n \frac{x}{1!} + \rho^{2n} \frac{x^2}{2!} + \cdots \right) = \rho^n e^{\rho^n x},
\end{aligned} \qquad (7)$$

где, при получении (7) учитывали, что $\frac{x^{-1}}{(-1)!} = \delta(x) = 0$ из за $x \geq x_0 > 0$.

Теперь интегрируем равенство (7) с дробным порядком $\frac{s+1}{n}$

$$I_{x_0}^{\frac{s+1}{n}} D_{x_0}^{\frac{s+1}{n}} J_s(x,\rho^n) = \rho^n I_{x_0}^{\frac{s+1}{n}} e^{\rho^n x},$$

т.е., по определению дробного интеграла Римана-Лиувилля [1], получаем

$$J_s(x,\rho^n) = \rho^n \int_{x_0}^{x} \frac{(x-t)^{\frac{s+1}{n}-1}}{\left(\frac{s+1}{n}-1\right)!} e^{\rho^n t} dt + \tilde{J}_s(x,\rho^n), \qquad (8)$$

где

$$D_{x_0}^{\frac{s+1}{n}} \tilde{J}_s(x, \rho^n) = 0. \tag{9}$$

Легко видеть, что из (9) можно определить $\tilde{J}_s(x, \rho^n)$ в следующем виде:

$$\tilde{J}_s(x, \rho^n) = e^{\rho^n x_0} \frac{x^{\frac{s+1}{n}-1}}{\left(\frac{s+1}{n}-1\right)!}. \tag{10}$$

## 4. Вычисление интеграла из (8)

Учитывая (10) в (8) имеем

$$J_s(x, \rho^n) = \rho^n \int_{x_0}^{x} \frac{(x-t)^{\frac{s+1}{n}-1}}{\left(\frac{s+1}{n}-1\right)!} e^{\rho^n t} dt + e^{\rho^n x_0} \frac{x^{\frac{s+1}{n}-1}}{\left(\frac{s+1}{n}-1\right)!}.$$

Произведем замену

$$J_s(x, \rho^n) = -\rho^n \int_{x-x_0}^{0} \frac{\xi^{\frac{s+1}{n}-1}}{\left(\frac{s+1}{n}-1\right)!} e^{\rho^n(x-\xi)} d\xi + e^{\rho^n x_0} \frac{x^{\frac{s+1}{n}-1}}{\left(\frac{s+1}{n}-1\right)!} =$$

$$= \rho^n e^{\rho^n x} \int_{0}^{x-x_0} \frac{\xi^{\frac{s+1}{n}-1}}{\left(\frac{s+1}{n}-1\right)!} \sum_{m=0}^{\infty} \frac{(-1)^m}{m!} \left(\rho^n \xi\right)^m d\xi + e^{\rho^n x_0} \frac{x^{\frac{s+1}{n}-1}}{\left(\frac{s+1}{n}-1\right)!} =$$

$$= \rho^n e^{\rho^n x} \frac{1}{\left(\frac{s+1}{n}-1\right)!} \sum_{m=0}^{\infty} \frac{(-1)^m}{m!} \rho^{nm} \int_{0}^{x-x_0} \xi^{\frac{s+1}{n}-1+m} d\xi + e^{\rho^n x_0} \frac{x^{\frac{s+1}{n}-1}}{\left(\frac{s+1}{n}-1\right)!} =$$

$$= \rho^n e^{\rho^n x} \frac{1}{\left(\frac{s+1}{n}-1\right)!} \sum_{m=0}^{\infty} \frac{(-1)^m}{m!} \rho^{nm} \left.\frac{\xi^{\frac{s+1}{n}+m}}{\frac{s+1}{n}+m}\right|_{\xi=0}^{x-x_0} + e^{\rho^n x_0} \frac{x^{\frac{s+1}{n}-1}}{\left(\frac{s+1}{n}-1\right)!} =$$

$$= \rho^n e^{\rho^n x} \frac{1}{\left(\frac{s+1}{n}-1\right)!} \sum_{m=0}^{\infty} \frac{(-1)^m}{m!} \rho^{nm} \frac{(x-x_0)^{\frac{s+1}{n}+m}}{\frac{s+1}{n}+m} + e^{\rho^n x_0} \frac{x^{\frac{s+1}{n}-1}}{\left(\frac{s+1}{n}-1\right)!}. \tag{11}$$

Подставляя (11) в (5) для $h_{\frac{1}{n}}(x, \rho)$ имеем

$$h_{\frac{1}{n}}(x,\rho) = \sum_{s=0}^{n-1} \rho^s J_s(x,\rho^n) = \sum_{s=0}^{n-1} \rho^{s+n} e^{\rho^n x} \frac{1}{\left(\frac{s+1}{n}-1\right)!} \sum_{k=0}^{\infty} \frac{(-1)^k}{k!} \frac{(x-x_0)^{\frac{s+1}{n}+k}}{\frac{s+1}{n}+k} \rho^{nk} +$$

$$+ \sum_{s=0}^{n-1} \rho^{s+n} e^{\rho^n x_0} \frac{x^{\frac{s+1}{n}-1}}{\left(\frac{s+1}{n}-1\right)!}. \qquad (12)$$

### 5. Замена параметра

Теперь в (12) заменим параметр $\rho$ на новый параметра $\lambda^{\frac{1}{m}}$ т.е., $\rho = \lambda^{\frac{1}{m}}$. Тогда имеем

$$h_{\frac{1}{n}}(x,\rho) = h_{\frac{1}{n}}(x,\lambda^{\frac{1}{m}}). \qquad (13)$$

Пусть $\alpha = \frac{m}{n}$, где $m$ и $n$ натуральные числа, и при $1 < m < n$. При $n = 1$ то и $m = 1$.

Вычислим следующую дробную производную:

$$D^{\alpha} h_{\frac{1}{n}}(x,\lambda^{\frac{1}{m}}) = D^{\frac{m}{n}} h_{\frac{1}{n}}(x,\lambda^{\frac{1}{m}}) = \lambda h_{\frac{1}{n}}(x,\lambda^{\frac{1}{m}}), \qquad (14)$$

Отметим, что (14) вытекает из (1) и (2).

Учитывая замену $\rho = \lambda^{\frac{1}{m}}$ в (12) имеем

$$h_{\frac{1}{n}}(x,\lambda^{\frac{1}{m}}) = \sum_{s=0}^{n-1} \lambda^{\frac{s+n}{m}} e^{\lambda^{\frac{n}{m}} x} \frac{1}{\left(\frac{s+1}{n}-1\right)!} \sum_{k=0}^{\infty} \frac{(-1)^k}{k!} \lambda^{\frac{nk}{m}} \frac{(x-x_0)^{\frac{s+1}{n}+k}}{\frac{s+1}{n}+k} +$$

$$+ \sum_{s=0}^{n-1} \lambda^{\frac{s}{m}} e^{\lambda^{\frac{n}{m}} x_0} \frac{x^{\frac{s+1}{n}-1}}{\left(\frac{s+1}{n}-1\right)!}. \qquad (15)$$

Таким образом, функция Миттаг-Леффлера со сдвигом (1) через экспоненциальную функцию представлена в виде (15). Для проверки результата (15) рассмотрим случай, $\alpha = 1$, т. е., $n = 1$, $m = 1$.

Тогда из (15) имеем:

$$h_{\frac{1}{n}}(x,\lambda) = \lambda e^{\lambda x}\sum_{k=0}^{\infty}\frac{(-1)^k}{k!}\lambda^k\frac{(x-x_0)^{1+k}}{1+k}+e^{\lambda x_0} = -e^{\lambda x}\sum_{k=0}^{\infty}\frac{(-1)^{k+1}}{(k+1)!}\lambda^{k+1}(x-x_0)^{k+1}+e^{\lambda x_0} =$$

$$= -e^{\lambda x}\sum_{p=1}^{\infty}\frac{(-1)^p}{(p)!}\lambda^p(x-x_0)^p+e^{\lambda x_0} = -e^{\lambda x}\left(e^{-\lambda(x-x_0)}-1\right)+e^{\lambda x_0} = e^{\lambda x}.$$

Таким образом, при $\alpha = 1$, функция Миттаг-Леффлера с сдвигом превращается к экспоненциальную функцию Эйлера $e^{\lambda x}$.

**6. Заключение:**

В работе впервые показывается что, функцию Миттаг-Леффлера со сдвигом можно представить через экспоненциальную функцию. В частном случае, когда порядок дробной производной равняется единице, то функция Миттаг-Леффлера со сдвигом превращается к экспоненциальную функцию Эйлера $e^{\lambda x}$.

# Representation of the Mittag-Leffler function through the exponential functions in the case of rational derivatives

**Fikret A. Aliev, N.A. Aliev, N.A. Safarova**

**Abstract.** In this paper the Mittag-Leffler function is given through the exponential functions for any rational derivatives of $\frac{m}{n}$ order, where $m < n$, $n > 1$ are natural irreducible numbers (if $n = 1$ then $m$ is also equal to unity). Unlike the previous papers the given formulas do not contain integrals.

**Keywords:** Mittag-Leffler function, exponential function, rational number.